\documentclass[12pt]{iopart}

\usepackage{graphicx,epstopdf}
\usepackage{amssymb}

\begin{document}

\title[A comparative study]{A comparative study of two globally convergent numerical methods for acoustic tomography}

\author{M V Klibanov$^1$ and A Timonov$^2$} 

\address{$^1$  Department of Mathematics, University of North Carolina at Charlotte, Charlotte, NC, USA}
\ead{mklibanv@math.ucf.edu}

\address{$^2$ Steklov Mathematical Institute of the Russian Academy of Sciences (S.Petersburg branch), S.Petersburg, Russia \\ University of South Carolina Upstate, Spartanburg, SC, USA}
\ead{altim@pdmi.ras.ru}

\vspace{10pt}
\begin{indented}
\item[]December 2022
\end{indented}

\begin{abstract}
The comparative study of two globally convergent numerical methods for acoustic tomography is carried out in two dimensions. These are the boundary control method and the quasi-reversibility method. The novelty is that in the latter a nonlinear inverse problem is reduced to a family of the linear integral equations of the first kind via the Lavrentiev approach and this reduction is used within the quasi-reversibility method. The analysis of its stability is carried out via Carleman estimates.  The computational effectiveness of these methods is tested in the numerical experiments with the smooth and discontinuous coefficients to be recovered from the tomographic data.    
\end{abstract}

%
%
%
%
%

\section{Introduction}
\label{intro}

In this paper, we study the computational effectiveness of two methods for the approximate solution of coefficient inverse problems as applied to acoustic tomography. These are the Boundary Control Method(BCM), originally proposed in \cite
{bel,bel-kyr}, and the regularized Quasi-Reversibility Method(QRM) based on the reduction of a nonlinear inverse problem for an acoustic wave equation to a linear integral equation \cite{klib, kokur, lavrmm}. Such a choice of methods is motivated by the fact that both of them transform the non-linear inverse problems to the linear ones. In addition, these methods provide, at least theoretically, the global convergence of approximate solutions to exact ones. In this paper, the acoustic tomography is interpreted as a research field for constructing and studying the numerical techniques for recovering the spatial distributions of some physical quantities, such as the speed of waves, the mass density, etc., from signals acquired by an acoustic transducers situated in a closed curve around an inhomogeneity to be investigated. 

The QRM was pioneered by Lattes and Lions in \cite{LL}. We prove the convergence of our version of the QRM using a Carleman
estimate. The first works where the convergence rates of the QRM were established are \cite{KS,KM,KR}. This was done via Carleman estimates. We also refer to the works \cite{B6}-\cite{ClK}, \cite{Ksurvey}, where Carleman
estimates were used to establish convergence rates of various versions of the QRM. The works \cite{QRMsurvey} and \cite[Chapter 4]{klib_b} contain surveys of the corresponding results.

In this paper, we use the QRM to solve the so-called Lavrent'ev integral equation of the first kind. This equation was first derived in the elegant work of Lavrent'ev in the 3-D case \cite{lavrmm}. One can also consider this equation as the 3-D analogue of the well known Gelfand-Levitan linear integral equation \cite{GL}, which solves a coefficient inverse problem for the 1-D wave equation The 2-D analogue of the Lavrent'ev equation was derived in \cite{kokur}. The first uniqueness
theorem for the Lavrent'ev equation in the non-overdetermined 3-D case was established by Kokurin \cite{kokurin}. The numerical method in the 3-D case, which is analogous to the one of this paper, was established in \cite{klib}.
However, our convergence analysis presented in this paper is different from the one indicated in \cite{klib}. The same is true for the numerical implementations. Also, see \cite{Bak1,Bak2} for two other numerical approaches to the 3-D Lavrent'ev equation.

In this paper, we work in $\mathbb{R}^2$. To model the propagation and scattering of acoustic waves, we utilize the simplest wave equation 
\begin{equation}
q(x)u_{tt}-\nabla ^{2}u=0~\mbox{in}~\Omega \times (0,\infty )  \label{simpl}
\end{equation}
where $\Omega \subset \mathbb{R}^{2}$ and $c=q^{-1/2}$ is the wave speed. In the quantitative acoustic imaging, it is required to recover, at least approximately, the distribution of the values of the wave speed $c$ (or the
coefficient $q$) in a bounded domain $\Omega $ from the measurements of the wave field $u(x,t)$ either on $\partial \Omega $ or on some subset of $\mathbb{R}^{2}\setminus\overline{\Omega }$. In this paper, all functions are supposed to be real valued. 

The paper is organized as follows. In the section~\ref{prelim} we formulate two inverse problems, which we consider. In the section~\ref{quant} we formulate two approaches to those inverse problems. A significant part of the section~\ref{stab} is devoted to the convergence analysis of our version of the QRM. A Carleman estimate for the Laplace operator is the key tool of the convergence analysis. The section~\ref{numer} is devoted to numerical studies. We conclude the paper in the section~\ref{concl}.

\section{Formulations of the inverse problems}
\label{prelim}

Let $\Omega\subset\mathbb{R}^{2}$ be a bounded and simply connected domain with the smooth boundary $\partial\Omega$, and $T>0$be the final time. As a first mathematical model of quantitative acoustic imaging, consider the initial boundary value problem
\begin{eqnarray}
& & q(x)u_{tt} - \Delta{u} = 0~\mbox{in}~\Omega\times(0, T],
\label{e_1} \\
& & u(x,0) = u_{t}(x,0) = 0~\mbox{in}~\overline{\Omega}\times\{0\}, \label{e_2} \\
& & \nabla{u}\cdot\nu = f~\mbox{on}~\partial\Omega\times[0,T], \label{e_3}
\end{eqnarray}
which has been extensively uzed in the mathematics literature to describe a dynamical system with the Neumann control. In this model, $f\in L^{2}(\partial\Omega\times[0,T])$ is the given control, $\nu$ is the outward unit normal vector on $\partial\Omega$, and the variable coefficient $q>0$ is supposed to be sufficiently smooth and bounded away from zero and infinity. Here and below, all functions in this model are supposed to be real valued. Denote $u^{f}$ the solution of (\ref{e_1})-(\ref{e_3}) that corresponds to $f$. If $q$ represents the mass density then the quantity $c = q^{-1/2}$ is the speed of waves. We choose the final time $T$ so that $T > sup\{d_{q}(x,\partial\Omega):x\in\Omega\}$, where $d_{q}$ is the distance in the metric $q^{1/2}|dx|$. Obviously, $\overline{\Omega} = \Omega^{T}$, where $\Omega^{T} = \{x\in\overline{\Omega}: d_{q}(x,\partial\Omega)\le T\}$. Under these conditions, the closure $\overline\Omega$ is spanned by the geodesic lines of $u^{f}$. Note that if 
$f\in C^{\infty}(\partial\Omega\times[0,T])$ and ${\it supp}f\subset(\partial\Omega\times(0,2T])$ then the solution $u^{f}$ is classical at $t=2T$. Based on the first model, we pose the following inverse problem. 

\noindent
{\bf Inverse Problem 1.} For the fixed $T > sup\{d_{q}(x,\partial\Omega):x\in\Omega\}$ given the set $\{(f,u^{f}_{|_{\partial\Omega\times(0,2T)}}): f\in L^{2}( \partial\Omega\times(0,2T) )\}$, determine $q$ in $\Omega^{T}$. 

\noindent
Note that since this set generates the Neumann-to-Dirichlet map $\Lambda^{2T}_{\partial\Omega}: f\rightarrow u^{f}$, the Inverse Problem 1 can also be stated as: Given $\Lambda^{2T}_{\partial\Omega}$, determine $q$ in $\Omega^{T}$.

Along with (\ref{e_1})-(\ref{e_3}), we also model the propagation and scattering of acoustic waves by the Cauchy problem
\begin{eqnarray}
& & q(x)u_{tt} - \Delta{u} = 0~\mbox{in}~\mathbb{R}^{2}\times(0, \infty),
\label{c_1} \\
& & u(x,0) =0,~~ u_{t}(x,0) = \delta(x-x_{0}),~x,x_{0}\in\mathbb{R}^{2}, \label{c_2} 
\end{eqnarray}
where $\delta(\cdot)$ is the Dirac delta function. 
Assume the coefficient $q$ is sufficiently smooth in $\mathbb{R}^{2}$. Besides, $q\ge 1$ in $\mathbb{R}^{2}$ and $q=1$ in $\mathbb{R}^{2}\setminus\Omega$. 
Note that this problem is equivalent to 
\begin{eqnarray}
& & q(x)u_{tt} - \Delta{u} = \delta(x-x_{0})\delta(t)~\mbox{in}~\mathbb{R}^{2}\times(0, \infty),
\label{c_11} \\
& & u(x,t)\equiv 0~\mbox{for}~t<0. \label{c_21} 
\end{eqnarray}
The solvability and stability results for this problem can be found in the mathematics literature (see,e.g., \cite{isakov, Rom}). This model gives rise to 

\noindent
{\bf Inverse Problem 2.} Without loss of generality, suppose that the origin is placed inside $\Omega$. Let $C_{R}$ be a circle with the center at the origin and with the radius $R$, such that 
%
\begin{equation}
G_{R}=\left\{ \left\vert x\right\vert <R\right\} ,\overline{\Omega }\subset
G_{R},C_{R}=\partial G_{R}.  \label{2.1}
\end{equation}
%
Let $u$ be the solution of (\ref{c_1})-(\ref{c_2}). Given the functions
\begin{eqnarray}
& & p_{0}(x,x_{0},t) = u(x,x_{0},t), \label{dat_1} \\
& & p_{1}(x,x_{0},t) = \nabla{u}(x,x_{0},t)\cdot\nu(x), \label{dat_2}  
\end{eqnarray}
where $(x,x_{0})\in C_{R}, x\neq x_{0}$, determine $q$ in $\Omega$, assuming that $q =1$ in 
$\mathbb{R}^{2}\diagdown \Omega$.

\section{Two approaches to the approximate solutions of the inverse problems} 
\label{quant}

\subsection{The Pestov's approach to the BCM}
\label{pest}

According to \cite{pest}, the Pestov's approach to the BCM is based on the approximate $H^{1}$-controllability of 
(\ref{e_1})-(\ref{e_3}). Denote $H_{q}$ the Hilbert space endowed with the norm
\[
\|w\|_{H_{q}} = \left( \int_{\Omega}qw^2 + \int_{\Omega}|\nabla{w}|^2 \right)^{1/2}.
\]
Let $u^{f}$ be a solution to (\ref{e_1})-(\ref{e_3}) that corresponds to the control $f$. Define a set
\[
{\cal U}^{T} = \left\{u^{f}(\cdot,T): f\in C^{\infty}(\partial\Omega\times[0,2T]), {\it supp}f\subset\partial\Omega\times[0,2T] \right\}.
\]
Let $T_{s} = \sup\{d_{q}(x,\partial\Omega): x\in\Omega \}$. For $T>T_{s}$ the reachable set of waves ${\cal U}^{T}$ is dense in $H_{q}$. That is, the following theorem takes place (see \cite{pest}, p.706).

{\bf Theorem 1.} {\it The orthogonal complement to ${\cal U}^{T}$ in $H_{q}$ is $\{0\}$.} 

This result allows for deriving an approximate relation 
\begin{equation}
\int_{\Omega}q(x)h_{1}(x)h_{2}(x)dx \approx [f_{h_{1}},f_{h_{2}}],
\label{intg}
\end{equation}
where $h_{1}, h_{2}$ are arbitrary harmonic functions and $[f_{h_{1}},f_{h_{2}}]$ is the bilinear form generated by these functions, which is defined as 
\[
[f_{h_{1}},f_{h_{2}}] = \int_{\Omega}q(x)u^{f_{h_{1}}}(x,T)u^{f_{h_{2}}}(x,T)dx.
\]
%
Since the linear span of all products $h_{1}h_{2}$ is dense in $L^{2}(\Omega)$, (\ref{intg}) is used for reducing the problem of obtaining an approximate $q$ to the solution of a system of linear algebraic equations. 
In theory, for any harmonic function $h\in C^{1}(\overline{\Omega})$ one can find the control $f\in L^{2}(\partial\Omega)$, such that $\|u^{f_{h}}(x,T) - h\|_{H^{1}(\Omega)}\le C\epsilon$, where $\epsilon$ is a small real number and $C=const>0$, which does not depend on $f$ and $h$. A numerical procedure for obtaining $h$ and $f$ is described in \cite{pest}.  

\subsection{Applying the method of quasi-reversibility}
\label{qrm}

Let the integrable function $\varphi(t)$, such that $\varphi=0$ for $t<0$, be of the exponential order, i.e., 
$|\varphi(t)|\le Ce^{p_{0}t}, C=const>0$ as $t\rightarrow\infty$. Then its Laplace transform is defined by
\begin{equation}
\tilde{\varphi}(p) = \int_{0}^{\infty}e^{-pt}\varphi(t)dt,~p > p_{0}, 
\label{lap}
\end{equation}
Here and below, we limit ourselves by the real values of $p$. Let $u$ be the solution of (\ref{c_1})-(\ref{c_2})(or (\ref{c_11})-(\ref{c_21})), and together with its derivatives $u_{tt}, u_{x_{1}x_{1}}, u_{x_{2}x_{2}}$ it has the Laplace transforms. According to the general Lavrentiev's approach \cite{lavrmm} and its 2D-implementation in \cite{kokur}, one may apply the Laplace transform to (\ref{c_1})-(\ref{c_2})(or (\ref{c_11})-(\ref{c_21})), introduce a new variable  
\begin{equation}
v(x,x_{0};p) = \frac{\tilde{u}(x,x_{0};p)}{\ln{p}},
\label{v_f}
\end{equation}
and reduce the problem (\ref{c_1})-(\ref{c_2})(or (\ref{c_11})-(\ref{c_21})) to a linear integral equation of the first kind. In the two dimensions this reduction was carried out in \cite{kokur}. For the reader's convenience, we summarize those results subject to our conditions. 
 
{\bf Theorem 2.}~{\it Suppose that the solution $u$ of (\ref{c_1})-(\ref{c_2})(or (\ref{c_11})-(\ref{c_21}) is such that for a certain integer $m\ge 0$ the function $(1+|x|^2 + t^2)^{-m}u(x,t)\in L^{1}(\mathbb{R}^3)$. Then $\xi$ satisfies the integral equation
\begin{equation}
\int_{\Omega}\ln{|x-x^{'}|}\ln{|x^{'}-x_{0}|}\xi(x^{'})dx^{'} = \psi(x,x_{0}),~x,~x_{0}\in C_{R},
\label{lavren}
\end{equation}
where
\begin{eqnarray}
& & \psi(x,x_{0}) = \lim_{p\rightarrow 0}(\ln{p}+\gamma)^{2}\left(
h(x,x_{0};p)-H_{0}-\frac{H_{1}(x,x_{0})}{\ln{p}+\gamma} \right), \label{rel_1} \\
& & H_{1}(x,x_{0}) = \lim_{p\rightarrow 0}(\ln{p}+\gamma)\left( h(x,x_{0};p)-H_{0} \right), \label{rel_2} \\
& & H_{0} = \lim_{p\rightarrow 0}h(x,x_{0};p), \label{rel_3} \\
& & h(x,x_{0};p) = 4\pi^{2}\cdot\frac{ v(x,x_{0};p) - g_{0}(x,x_{0};p) }{ (1+\gamma(\ln{p})^{-1})^{2}p^{2}\ln{p} }, \label{rel_4} \\
& & g_{0}(x,x_{0};p) =  -\frac{1}{\pi}\left( 1+\frac{\gamma}{\ln{p}} + \frac{1}{\ln{p}}\ln{|x-x_{0}|} +  |x-x_{0}|^{2}p^2 \right) +   \nonumber \\ 
& & + \left( \frac{p^{2}|x-x_{0}|^2}{\ln{p}}\ln{ |x-x_{0}| } + (\gamma -1)\frac{p^{2}|x-x_{0}|^2}{\ln{p}}\right). \label{rel_5}
\end{eqnarray}
}
Here, $\xi(x) = q(x) - 4$ and $p\in(0,\varepsilon), \varepsilon\in(0,e^{-\gamma})$, $\gamma = 0.5772...$ is the Euler constant.

Assume that $q\in C^{\beta}, \beta\in(0,1)$ and denote
\[
u(x,x_{0}) = \frac{1}{2\pi}\int_{\Omega}\ln{|x-x^{'}|}\ln{|x^{'}-x_{0}|}\xi(x^{'})dx^{'},~x\in\Omega,~x_{0}\in C_{R}.
\]
Applying the operator $\Delta_{x}$ to the both sides of (\ref{lavren}), we then obtain
\begin{equation}
\Delta_{x}u(x,x_{0}) = \xi(x)\ln{|x-x_{0}|},~x\in\Omega,~x_{0}\in C_{R}.
\label{eq_x}
\end{equation}
subject to the boundary conditions on $C_{R}, x\neq x_{0}$
\begin{eqnarray}
& & u(x,x_{0}) = \psi(x,x_{0}), \label{ch_1} \\
& & \nabla{u}\cdot\nu = \psi_{1}(x,x_{0}). \label{ch_2}
\end{eqnarray}
%
Note that in order to find $\psi_1$ on $C_{R}$, it is sufficient to solve the exterior boundary value problem, i.e., to solve 
(\ref{eq_x}) with $\xi=0$ in $\mathbb{R}^{2}\setminus\overline{\Omega}$ subject to (\ref{ch_1}) and $u\rightarrow 0$ as $|x|\rightarrow\infty$.
Then the Inverse Problem 2 can be reformulated as follows. Assume $u(x,x_{0})$ satisfies (\ref{eq_x}). Given $\psi(x,x_{0}),\psi_{1}(x,x_{0})$ on $C_{R}$, find $\xi$ in $\Omega$.    

Next, we introduce a new variable
\[
v(x,x_{0}) = \frac{u(x,x_{0})}{\ln{|x-x_{0}|}},~x\in\Omega,~x_{0}\in C_{R}.
\]
Substituting it in (\ref{eq_x}) and resolve the resulting equation for $\xi$, we obtain
\begin{equation}
\xi(x) = \Delta_{x}v - 2\nabla_{x}(v\ln{|x-x_{0}|})\nabla_{x}\left(\frac{1}{\ln{|x-x_{0}|}}\right) - v\left[\ln{|x-x_{0}|}\nabla_{x}\left(\frac{1}{\ln{|x-x_{0}|}}\right)\right].
\label{diff}
\end{equation}
Since $\xi$ does not depend on $x_{0}$, we differentiate (\ref{diff}) with respect to $x_{0}$ and obtain the partial differential equation of the third order that does not contain $\xi$
\begin{eqnarray}
& & \Delta_{x}v_{x_{0}} - 2\frac{\partial}{\partial x_{0}}\left\{\nabla_{x}(v\ln{|x-x_{0}|})\nabla_{x}\left(\frac{1}{\ln{|x-x_{0}|}}\right)\right\} - \nonumber \\ 
& & -\frac{\partial}{\partial x_{0}}\left\{v\left[\ln{|x-x_{0}|}\nabla_{x}\left(\frac{1}{\ln{|x-x_{0}|}}\right)\right]\right\} = 0 \label{eq_main}
\end{eqnarray}
subject to the boundary conditions
\begin{eqnarray}
& & v(x,x_{0}) = \frac{\psi(x,x_{0})}{\ln{|x-x_{0}|}} = s_{0}(x,x_{0}), \label{bc_1} \\
& & \nabla{v(x,x_{0})}\cdot\nu(x) = \frac{\psi_{1}(x,x_{0})}{\ln{|x-x_{0}|}} - \frac{\psi(x,x_{0})}{\ln{|x-x_{0}|}}\frac{\partial}{\partial\nu}(\ln{|x-x_{0}|}). \label{bc_2}
\end{eqnarray}
It is easy to see that Inverse Problem 2 is actually reduced to the problem of extending the field $v(x,x_{0})$ from the boundary $C_{R}$ into $\Omega$.  

Now, we apply the quasi-reversibility method to the approximate solution of (\ref{eq_main})-(\ref{bc_2}). To do this, consider a special orthonormal basis $\{\Psi_{k}(x_{0})\}_{k=0}^{\infty}$ (see \cite{klib_b}) in $L^{2}(C_{R})$ and approximate the function $v(x,x_{0})$ by 
\begin{equation}
v(x,x_{0})\approx\sum_{k=0}^{N-1}v_{k}(x)\Psi_{k}(x_{0}), \label{decomp}
\end{equation}
where
\[
v_{k}(x) = \int_{C_{R}}v(x,x_{0})\Psi_{k}(x_{0})dx_{0},~k=0,...,N-1.
\]
Also, we approximate the functions
\begin{eqnarray}
& & s_{0}(x,x_{0})\approx\sum_{k=1}^{N-1}s_{0k}(x)\Psi_{k}(x_{0}), \label{bc1_a} \\
& & s_{1}(x,x_{0})\approx\sum_{k=1}^{N-1}s_{1k}(x)\Psi_{k}(x_{0}), \label{bc2_a}
\end{eqnarray}
where
\begin{eqnarray*}
& & s_{0k}(x) = \int_{C_{R}}s_{0}(x,x_{0})\Psi_{k}(x_{0})dx_{0}, \\
& & s_{1k}(x) = \int_{C_{R}}s_{1}(x,x_{0})\Psi_{k}(x_{0})dx_{0},~k=0,...,N-1.
\end{eqnarray*}
We substitute (\ref{decomp}) in (\ref{eq_main}) and multiply subsequently both sides of the resulting equation by 
$\Psi_{k}(x_{0}), k=0,...,N-1$, and then integrate with respect to $x_{0}\in C_{R}, x\neq x_{0}$. As a result, we obtain a system of elliptic PDEs 
\begin{equation}
L(V)\equiv M_{N}\Delta{V} + K_{1}(x)V_{x_{1}} + K_{2}(x)V_{x_{2}} + K_{0}(x)V = 0,~x=(x_{1},x_{2}) \label{syst_V}
\end{equation}
subject to the boundary conditions
\begin{equation}
V=S_{0}(x),~~\frac{\partial V}{\partial\nu(x)}=S_{1}(x),~x\in C_{R}.
\label{bc_V}
\end{equation}
Here, $V = (v_{0},...,v_{N-1})^{t}(x), S_{0} = (s_{00},...,v_{0N-1})^{t}(x), S_{1} = (s_{10},...,s_{1N-1})^{t}(x)$ 
are the pseudo-vectors and $K_{j}, j=0,1,2$ are the $N\times N$-pseudo-matrices. The matrix $M_{N}$ is given by $M_{N} = (a_{kn})^{N-1}_{k,n=0}$,
\begin{equation}\label{matr}
a_{kn}=\left\{
\begin{array}{l}
1 \mbox{ if } k=n\\
0 \mbox{ if } k>n, 
\end{array}\right.
\end{equation}
By virtue of (\ref{matr}), ${\it det}M_{N} = 1$, i.e. there exists an inverse $M^{-1}_{N}$, though the condition number of $M_{N}$ increases rapidly as $N$ grows. Once the problem (\ref{syst_V})-(\ref{bc_V}) is approximately solved, an approximation of 
$\xi$ is determined from (\ref{diff}). In this connection, one may interpret the aforementioned procedure as a technique for obtaining an approximate single internal data for the Inverse Problem 2. 

Since the problem (\ref{syst_V})-(\ref{bc_V}) is over-determined, the Tikhonov's smoothing functional together with the Carleman estimates are utilized in the next section for the purpose of analysis.     

\section{The convergence analysis}
\label{stab}

To begin with the convergence analysis, consider the Tikhonov's smoothing functional with the
regularization parameter $\alpha \in \left( 0,1\right)$ 
\begin{equation}
J_{\alpha }(V)= \int_{G_{R}}(L(V))^{2}dx+\alpha \Vert V\Vert
_{H^{2}(G_{R})}^{2} \label{0}
\end{equation}
on the set 
\begin{equation}
B=\left\{ V\in H^{2}(G_{R}):V=S_{0},\frac{\partial }{\partial \nu }V=S_{1}~
\mbox{on}~C_{R}\right\} .  \label{0.1}
\end{equation}
Here and below, $V\in H^{2}(G_{R})$ means that each component of the vector
function $V$ belongs to the space $H^{2}(G_{R})$ and 
\[
\left\Vert V\right\Vert
_{H^{2}(G_{R})}^{2}=\sum\limits_{n=0}^{N-1}\left\Vert v_{n}\right\Vert
_{H^{2}(G_{R})}^{2}.
\]
We solve 

\noindent
\textbf{The minimization problem:} \emph{Minimize the functional }$J_{\alpha }(V)
$\emph{\ in (\ref{0}) on the set }$B$\emph{.} 

\noindent
\textbf{Remark 3.1.} \emph{Suppose that we have found a minimizer }$V_{\min
}\in B$\emph{\ of functional (\ref{0}). Then to find an approximation to the
unknown function }$\xi \left( x\right) ,$\emph{\ we substitute in the right
hand side of (\ref{decomp}) components of the vector function }$V_{\min }.$%
\emph{\ Next, we substitute the resulting left hand side of (\ref{decomp})
in the right hand side of (\ref{diff}). Since the right hand side of (\ref%
{diff}) depends on }$x_{0},$\emph{\ then one needs to integrate it with
respect to }$x_{0}\in C_{R}$\emph{\ and then divide by }$2\pi R$\emph{,
which is the length of the circle }$C_{R}$\emph{.Then the left hand side of
the resulting formula will be the needed approximation for the function }$%
\xi \left( x\right) $\emph{.}

\subsection{Existence and uniqueness of the minimizer}

\textbf{Theorem 3.}\emph{\ Suppose that there exists a vector function }$W\in H^{2}(G_{R})$ such that 
\begin{equation}
W=S_{0},\frac{\partial }{\partial \nu }W=S_{1}~\mbox{on}~C_{R}.  \label{1}
\end{equation}
\emph{Then there exists unique minimizer }$V_{\min }\in B$\emph{\ of the
functional }$J_{\alpha }(V)$\emph{\ on the set }$B$\emph{\ and the following
estimate holds:}
\begin{equation}
\left\Vert V_{\min }\right\Vert _{H^{2}(G_{R})}\leq \frac{C}{\sqrt{\alpha }}%
\left\Vert W\right\Vert _{H^{2}(G_{R})},  \label{2}
\end{equation}
\emph{where the number }$C=C\left( L,G_{R}\right) >0$\emph{\ depends only on
listed parameters.}

\textbf{Proof.} For each vector function $V\in B$ consider the vector
function $\widetilde{V}$ defined as
\begin{equation}
\widetilde{V}=V-W.  \label{3}
\end{equation}
Define 
\begin{equation}
\emph{\ }H_{0}^{2}\left( \Omega \right) =\left\{ P=\left(
p_{0},...,p_{N-1}\right) ^{T}\in H^{2}\left( G_{R}\right) :P=\frac{\partial 
}{\partial \nu }P=0~\mbox{on}~C_{R}\right\} .  \label{3.0}
\end{equation}
Consider the functional $I_{\alpha }:H_{0}^{2}\left( G_{R}\right)
\rightarrow \mathbb{R}$ defined as
\begin{equation}
I_{\alpha }\left( P\right) =J_{\alpha }\left( P+W\right) ,\forall P\in
H_{0}^{2}\left( G_{R}\right) .  \label{4}
\end{equation}
Suppose that there exists a minimizer $V_{\min }\in B$\emph{\ }of the
functional $J_{\alpha }(V).$ Then by (\ref{3}) and (\ref{4})
\[
I_{\alpha }\left( \widetilde{V}_{\min }\right) =J_{\alpha }\left( \widetilde{%
V}_{\min }+W\right) =J_{\alpha }\left( V_{\min }\right) \leq J_{\alpha
}\left( \widetilde{V}+W\right) =I_{\alpha }\left( \widetilde{V}\right) ,
\]
\begin{equation}
\forall \widetilde{V}=V-W,\forall V\in B.  \label{5}
\end{equation}
Consider now an arbitrary $P\in H_{0}^{2}\left( G_{R}\right) .$ Then $P+W=Q\in B.$ Hence, by (\ref{3}) $P=\widetilde{Q}=Q-W.$ Hence, by (\ref{5})
\[
I_{\alpha }\left( \widetilde{V}_{\min }\right) \leq I_{\alpha }\left(
P\right) ,\forall P\in H_{0}^{2}\left( G_{R}\right) .
\]%
Therefore, we have proven that $\widetilde{V}_{\min }$ is a minimizer of
functional (\ref{4}) on the set $B$. Reversed is also obviously true: If a
vector function $Y_{\min }\in H_{0}^{2}\left( G_{R}\right) $ is a minimizer
of functional (\ref{4}) on the space $H_{0}^{2}\left( G_{R}\right) ,$ then $
Y+W\in B$ is a minimizer of the functional $J_{\alpha }(V)$ on the set $B$.
Therefore, we consider now the problem of the minimization of the functional 
$I_{\alpha }\left( P\right) $ on the whole space $H_{0}^{2}\left(
G_{R}\right) .$ 

By the variational principle and (\ref{0}) we should have%
\begin{equation}
\left( L\left( h\right) ,L\left( Y_{\min }\right) \right) +\alpha \left[
h,Y_{\min }\right] =-\left( L\left( h\right) ,L\left( W\right) \right)
-\alpha \left[ h,W\right] ,\forall h\in H_{0}^{2}\left( G_{R}\right) ,
\label{6}
\end{equation}
where $\left( ,\right) $ and $\left[ ,\right] $ are scalar products in
spaces $L^{2}\left( G_{R}\right) $ and $H^{2}\left( G_{R}\right) $
respectively. Define a new scalar product in $H^{2}\left( G_{R}\right) ,$%
\begin{equation}
\left\{ f,g\right\} =\left( Lf,Lg\right) +\alpha \left[ h,g\right] ,\forall
f,g\in H_{0}^{2}\left( G_{R}\right) .  \label{7}
\end{equation}
Then obviously
\begin{equation}
\alpha \left\Vert f\right\Vert _{H^{2}\left( G_{R}\right) }^{2}\leq \left\{
f,f\right\} \leq C\left\Vert f\right\Vert _{H^{2}\left( G_{R}\right) }^{2}.
\label{8}
\end{equation}
Hence, the scalar product (\ref{7}) generates the norm $\left\{ \cdot
\right\} $ in $H_{0}^{2}\left( G_{R}\right) ,$ which is equivalent to the
original norm of the space $H^{2}\left( G_{R}\right) .$

Consider the right hand side of (\ref{6}). By Cauchy-Schwarz inequality 
\[
\left\vert -\left( L\left( h\right) ,L\left( W\right) \right) -\alpha \left[
h,W\right] \right\vert \leq C\left\Vert W\right\Vert _{H^{2}\left(
G_{R}\right) }\cdot \left\{ h\right\} ,\forall h\in H_{0}^{2}\left(
G_{R}\right) .
\]
Hence, the right hand side of (\ref{6}) can be considered as a bounded
linear functional defined on the space $H_{0}^{2}\left( G_{R}\right)$
\begin{equation}
l_{W}\left( h\right) =-\left( L\left( h\right) ,L\left( W\right) \right)
-\alpha \left[ h,W\right] ,\left\Vert l_{W}\right\Vert \leq C\left\Vert
W\right\Vert _{H^{2}\left( G_{R}\right) }.  \label{10}
\end{equation}
Hence, Riesz theorem (\ref{8}) and (\ref{10}) imply that there exists unique 
$U\in H_{0}^{2}\left( G_{R}\right) $ such that 
\begin{equation}
l_{W}\left( h\right) =\left\{ U,h\right\} ,\forall h\in H_{0}^{2}\left(
G_{R}\right) ,\left\Vert U\right\Vert _{H^{2}\left( G_{R}\right) }\leq
C\left\Vert W\right\Vert _{H^{2}\left( G_{R}\right) }.  \label{11}
\end{equation}
Hence, by (\ref{6}) and (\ref{7})
\[
\left\{ Y_{\min },h\right\} =\left\{ U,h\right\} ,\forall h\in
H_{0}^{2}\left( G_{R}\right) .
\]
Hence, the minimizer $Y_{\min }$ of the functional $I_{\alpha }\left(
P\right) $ on the space $H_{0}^{2}\left( G_{R}\right) $ indeed exists, is
unique and $Y_{\min }=U.$ Furthermore, by (\ref{8}) and (\ref{11}) 
\[
\left\Vert Y_{\min }\right\Vert _{H^{2}\left( G_{R}\right) }\leq \frac{C}{%
\sqrt{\alpha }}\left\Vert W\right\Vert _{H^{2}\left( G_{R}\right) }.
\]
Finally, setting $V_{\min }=Y_{\min }+W,$ we obtain the desired result,
including estimate (\ref{2}). $\square $

\subsection{Convergence of regularized solutions}

The minimizer of Theorem 3 is called the "regularized
solution" in the theory of regularization \cite{book_tgsy}. It is always
important to prove convergence of regularized solutions to the exact one as
long as the level of noise in the data tends to zero \cite{book_tgsy}. This
question is a more delicate one than the question of existence and
uniqueness of the regularized solution addressed in Theorem 3. Indeed, while
in Theorem 3 we have relied on the classic Riesz theorem, now we need to
bring in the tool of Carleman estimates.

Let the number $\mu \in \left( 0,R/2\right) .$ Using (\ref{2.1}), define 
\begin{equation}
G_{\mu ,R}=\left\{ x\in \mathbb{R}^{2}:\mu <\left\vert x\right\vert
<R\right\} .  \label{12}
\end{equation}
We assume that 
\begin{equation}
\Omega \subset G_{2\mu ,R}.  \label{13}
\end{equation}
Introduce polar coordinates 
\[
x_{1}=r\cos \varphi ,x_{2}=R\sin \varphi ,r\in \left( \mu ,R\right) ,\varphi
\in \left[ 0,2\pi \right) .
\]

\textbf{Theorem 4} (Carleman estimate \cite{KEIT}, \cite[Theorem 7.3.2]%
{klib_b}). \emph{There exists a number }$\lambda _{0}=\lambda _{0}\left(
G_{\mu ,R}\right) \geq 1$\emph{\ and a \ number }$C=C\left( G_{\mu
,R}\right) >0,$ \emph{both depending only on the domain }$G_{\mu ,R},$\emph{%
\ such that for all }$w\in H_{0}^{2}\left( G_{\mu ,R}\right) $\emph{\ and
for all }$\lambda \geq \lambda _{0}$\emph{\ the following Carleman estimate
holds:}
\[
\int\limits_{G_{\mu ,R}}\left( \Delta w\right) ^{2}e^{2\lambda r}dx\geq 
\frac{1}{2}\int\limits_{G_{\mu ,R}}\left( \Delta w\right) ^{2}e^{2\lambda
r}+C\lambda\int\limits_{G_{\mu ,R}}\left[ \left( \nabla w\right)
^{2}+\lambda ^{2}w^{2}\right] e^{2\lambda r}dx-
\]
\[
-C\lambda ^{3}e^{2\lambda \mu }\left\Vert w\right\Vert _{H^{2}\left( G_{\mu
,R}\right) }^{2}.
\]

The proof of Lemma 1 is not presented here since this lemma is well known.

\textbf{Lemma 1. }\emph{Let }$A$\emph{\ be an }$n\times n$\emph{\ matrix
such that the matrix }$A^{-1}$\emph{\ exists. Then there exists a number }$%
\gamma =\gamma \left( A\right) >0$\emph{\ depending only on }$A$\emph{\ such
that }$\left\Vert Ax\right\Vert \geq \gamma \left\Vert x\right\Vert ,\forall
x\in R^{n},$\emph{\ where norm is the euclidian norm.}

According to the theory of regularization, we assume now that there exists
exact solution $V^{\ast }\in H^{2}\left( G_{R}\right) $ of problem (\ref%
{syst_V})-(\ref{bc_V}) with the noiseless data $S_{0}^{\ast },S_{1}^{\ast },$
\begin{equation}
L\left( V^{\ast }\right) =0,  \label{130}
\end{equation}
\begin{equation}
V^{\ast }=S_{0}^{\ast },\frac{\partial }{\partial \nu }V^{\ast }=S_{1}^{\ast
}~\mbox{on}~C_{R},  \label{14}
\end{equation}
see (\ref{syst_V})-(\ref{bc_V}). Let the vector function $W\in H^{2}\left(
G_{R}\right) $ be the one satisfying (\ref{1}). Let $\delta \in \left(
0,1\right) $ be the level of noise in the data $S_{0},S_{1}$ in (\ref{0.1}).
We assume that there exists such an extention $W^{\ast }\in H^{2}\left(
G_{R}\right) $ of the noiseless data $S_{0}^{\ast },S_{1}^{\ast }$ in (\ref%
{14}) that 
\begin{equation}
W^{\ast }=S_{0}^{\ast },\frac{\partial }{\partial \nu }W^{\ast }=S_{1}^{\ast
}~\mbox{on}~C_{R},  \label{15}
\end{equation}%
\begin{equation}
\left\Vert W-W^{\ast }\right\Vert _{H^{2}\left( G_{R}\right) }<\delta .
\label{16}
\end{equation}

\textbf{Theorem 5.} \emph{Let the regularization parameter }$\alpha =\alpha
\left( \delta \right) =\delta ^{2}.$\emph{\ Suppose that (\ref{1}), (\ref{12}%
)-(\ref{16}) hold. Let }$V_{\min }\in B$\emph{\ be the of the functional }$%
J_{\alpha }(V)$\emph{\ on the set }$B,$\emph{\ which was found in Theorem 3.
There exists a sufficiently small number }$\delta _{0}=\delta _{0}\left(
L,G_{\mu ,R}\right) \in \left( 0,1\right) $\emph{\ and a number }$%
C_{1}=C_{1}\left( L,G_{\mu ,R}\right) >0,$\emph{\ both depending only on
listed parameters such that for all }$\delta \in \left( 0,\delta _{0}\right) 
$\emph{\ the following estimates hold}
\begin{equation}
\left\Vert \Delta V_{\min }-\Delta V^{\ast }\right\Vert _{L^{2}\left(
G_{2\mu ,R}\right) }+\left\Vert V_{\min }-V^{\ast }\right\Vert _{H^{1}\left(
G_{2\mu ,R}\right) }\leq   \label{160}
\end{equation}
\[
\leq C_{1}\left( 1+\left\Vert V^{\ast }\right\Vert _{H^{2}\left(
G_{R}\right) }\right) \delta ^{\mu /\left( 4R\right) },
\]
\begin{equation}
\left\Vert \xi _{\min }-\xi ^{\ast }\right\Vert _{L^{2}\left( \Omega \right)
}\leq C_{1}\left( 1+\left\Vert V^{\ast }\right\Vert _{H^{2}\left(
G_{R}\right) }\right) \delta ^{\mu /\left( 4R\right) },  \label{161}
\end{equation}
\emph{where functions }$\xi _{\min }\left( x\right) $\emph{\ and }$\xi
^{\ast }\left( x\right) $\emph{\ are reconstructed from vector functions }$%
V_{\min }$\emph{\ and }$V^{\ast }$\emph{\ respecively as in Remark 3.1.}

\textbf{Proof.} As in the proof of Theorem 3, let 
\begin{equation}
Y_{\min }=V_{\min }-W.  \label{17}
\end{equation}%
$.$ Also, let 
\begin{equation}
\widetilde{V}^{\ast }=V^{\ast }-W^{\ast }.  \label{18}
\end{equation}
Then by (\ref{1}), (\ref{14}), (\ref{15}), (\ref{17}) and (\ref{18}) both
vector functions $Y_{\min },\widetilde{V}^{\ast }\in H_{0}^{2}\left(
G_{R}\right) .$ \ Also, by (\ref{130})%
\begin{equation}
L\left( \widetilde{V}^{\ast }\right) =-L\left( W^{\ast }\right) .  \label{19}
\end{equation}

Using (\ref{6}) and (\ref{19}), we obtain
\[
\left( L\left( h\right) ,L\left( Y_{\min }\right) \right) +\alpha \left[
h,Y_{\min }\right] =
\]
\begin{equation}
=-\left( L\left( h\right) ,L\left( W\right) \right) -\alpha \left[ h,W\right]
,\forall h\in H_{0}^{2}\left( G_{R}\right) ,  \label{20}
\end{equation}
\[
\left( L\left( h\right) ,L\left( \widetilde{V}^{\ast }\right) \right)
+\alpha \left[ h,\widetilde{V}^{\ast }\right] =
\]
\begin{equation}
=-\left( L\left( h\right) ,L\left( W^{\ast }\right) \right) +\alpha \left[ h,%
\widetilde{V}^{\ast }\right] ,\forall h\in H_{0}^{2}\left( G_{R}\right) .
\label{21}
\end{equation}
Denote 
\begin{equation}
Q=Y_{\min }-\widetilde{V}^{\ast }\in H_{0}^{2}\left( G_{R}\right) .
\label{22}
\end{equation}
Subtract (\ref{21}) from (\ref{20}). Using (\ref{22}), we obtain
\[
\left( L\left( h\right) ,L\left( Q\right) \right) +\alpha \left[ h,Q\right] =
\]
\begin{equation}
=-\left( L\left( h\right) ,L\left( W-W^{\ast }\right) \right) -\alpha \left[
h,W+V^{\ast }\right] ,\forall h\in H_{0}^{2}\left( G_{R}\right) .  \label{23}
\end{equation}
Set in (\ref{23}) $h=Q$ and use Cauchy-Schwarz inequality and (\ref{16}). We
obtain%
\[
\left\Vert L\left( Q\right) \right\Vert _{L_{2}\left( G_{R}\right)
}^{2}+\alpha \left\Vert Q\right\Vert _{H^{2}\left( G_{R}\right) }^{2}\leq 
\frac{1}{2}\left\Vert L\left( Q\right) \right\Vert _{L_{2}\left(
G_{R}\right) }^{2}+C\delta ^{2}+\frac{\alpha }{2}\left\Vert Q\right\Vert
_{H^{2}\left( G_{R}\right) }^{2}+
\]%
\begin{equation}
+\frac{\alpha }{2}\left\Vert W+\widetilde{V}^{\ast }\right\Vert
_{H^{2}\left( G_{R}\right) }^{2}.  \label{24}
\end{equation}%
By triangle inequality, (\ref{16}) and (\ref{18})%
\begin{equation}
\left\Vert W+\widetilde{V}^{\ast }\right\Vert _{H^{2}\left( G_{R}\right)
}=\left\Vert W-W^{\ast }+V^{\ast }\right\Vert _{H^{2}\left( G_{R}\right)
}\leq \delta +\left\Vert V^{\ast }\right\Vert _{H^{2}\left( G_{R}\right) }.
\label{25}
\end{equation}%
Since $\alpha \in \left( 0,1\right) ,$ then (\ref{24}) and (\ref{25})\ imply%
\[
\left\Vert L\left( Q\right) \right\Vert _{L_{2}\left( G_{R}\right)
}^{2}+\alpha \left\Vert Q\right\Vert _{H^{2}\left( G_{R}\right) }^{2}\leq
C\delta ^{2}+2\alpha \left\Vert V^{\ast }\right\Vert _{H^{2}\left(
G_{R}\right) }^{2}.
\]%
Recalling that $\alpha =\delta ^{2},$ we obtain%
\begin{equation}
\left\Vert Q\right\Vert _{H^{2}\left( G_{R}\right) }^{2}\leq C\left(
1+\left\Vert V^{\ast }\right\Vert _{H^{2}\left( G_{R}\right) }^{2}\right) 
\label{26}
\end{equation}%
\begin{equation}
\int\limits_{G_{R}}\left[ L\left( Q\right) \right] ^{2}dx\leq C\delta
^{2}\left( 1+\left\Vert V^{\ast }\right\Vert _{H^{2}\left( G_{R}\right)
}^{2}\right) .  \label{27}
\end{equation}

What we have used in this proof so far were just considerations of the 
standard analysis of Hilbert spaces. Now, however, we are ready to use a
non-trivial \ argument of Carleman estimates. 

We have
\[
\int\limits_{G_{R}}\left[ L\left( Q\right) \right] ^{2}dx=\int
\limits_{G_{R}}\left[ L\left( Q\right) \right] ^{2}e^{2\lambda
r}e^{-2\lambda r}dx\geq e^{-2\lambda R}\int\limits_{G_{R}}\left[ L\left(
Q\right) \right] ^{2}e^{2\lambda r}e^{-2\lambda r}dx.
\]
Hence, (\ref{27}) implies
\begin{equation}
\int\limits_{G_{R}}\left[ L\left( Q\right) \right] ^{2}e^{2\lambda r}dx\leq
C\delta ^{2}e^{2\lambda R}\left( 1+\left\Vert V^{\ast }\right\Vert
_{H^{2}\left( G_{R}\right) }^{2}\right) .  \label{28}
\end{equation}
We now apply the Carleman estimate of Theorem 4 to the left hand side of (%
\ref{28}). This is possible since by (\ref{22}) $Q\in H_{0}^{2}\left(
G_{R}\right) .$ By  (\ref{syst_V}), Theorem 4, Lemma 1 and (\ref{28}) 
\[
C\delta ^{2}e^{2\lambda R}\left( 1+\left\Vert V^{\ast }\right\Vert
_{H^{2}\left( G_{R}\right) }^{2}\right) \geq \int\limits_{G_{\mu R}}\left[
L\left( Q\right) \right] ^{2}e^{2\lambda r}dx\geq 
\]
\[
\geq\int\limits_{G_{\mu R}}\left[ M_{N}\Delta Q\right] ^{2}e^{2\lambda
r}dx-C_{1}\int\limits_{G_{\mu R}}\left[ \left( \nabla Q\right) ^{2}+Q^{2}%
\right] e^{2\lambda r}dx\geq 
\]
\[
\geq \frac{1}{2}\int\limits_{G_{\mu R}}\left[ M_{N}\Delta Q\right]
^{2}e^{2\lambda r}dx+C_{1}\int\limits_{G_{\mu R}}\left[ \Delta Q\right]
^{2}e^{2\lambda r}dx-C_{1}\int\limits_{G_{\mu R}}\left[ \left( \nabla
Q\right) ^{2}+Q^{2}\right] e^{2\lambda r}dx\geq 
\]
\[
\geq C_{1}\int\limits_{G_{\mu R}}\left[ \Delta Q\right] ^{2}e^{2\lambda
r}dx+C_{1}\lambda\int\limits_{G_{\mu ,R}}\left[ \left( \nabla Q\right)
^{2}+\lambda ^{2}Q^{2}\right] e^{2\lambda r}dx-C\lambda ^{3}e^{2\lambda \mu
}\left\Vert Q\right\Vert _{H^{2}\left( G_{\mu ,R}\right) }^{2}-
\]
\[
-C_{1}\int\limits_{G_{\mu R}}\left[ \left( \nabla Q\right) ^{2}+Q^{2}\right]
e^{2\lambda r}dx.
\]%
Hence, there exists a sufficiently large $\lambda _{1}=\lambda _{1}\left(
L,M_{N},G_{\mu ,R}\right) \geq \lambda _{0}\left( G_{\mu ,R}\right) \geq 1$
such that for all $\lambda \geq \lambda _{1}$
\begin{equation}
\int\limits_{G_{\mu ,R}}\left[ \left( \Delta Q\right) ^{2}+\left( \nabla
Q\right) ^{2}+Q^{2}\right] e^{2\lambda r}dx\leq   \label{29}
\end{equation}%
\[
\leq C_{1}\delta ^{2}e^{2\lambda R}\left( 1+\left\Vert V^{\ast }\right\Vert
_{H^{2}\left( G_{R}\right) }^{2}\right) +C_{1}\lambda ^{3}e^{2\lambda \mu
}\left\Vert Q\right\Vert _{H^{2}\left( G_{\mu ,R}\right) }^{2}.
\]
Now, 
\[
\int\limits_{G_{\mu ,R}}\left[ \left( \Delta Q\right) ^{2}+\left( \nabla
Q\right) ^{2}+Q^{2}\right] e^{2\lambda r}dx\geq\int\limits_{G_{2\mu ,R}}%
\left[ \left( \Delta Q\right) ^{2}+\left( \nabla Q\right) ^{2}+Q^{2}\right]
e^{2\lambda r}dx\geq 
\]
\[
\geq e^{4\mu \lambda }\int\limits_{G_{2\mu ,R}}\left[ \left( \Delta
Q\right) ^{2}+\left( \nabla Q\right) ^{2}+Q^{2}\right] e^{2\lambda r}dx.
\]%
Hence, (\ref{29}) implies for sufficiently large $\lambda $
\[
\int\limits_{G_{2\mu ,R}}\left[ \left( \Delta Q\right) ^{2}+\left( \nabla
Q\right) ^{2}+Q^{2}\right] e^{2\lambda r}dx\leq C_{1}\delta ^{2}e^{2\lambda
R}\left( 1+\left\Vert V^{\ast }\right\Vert _{H^{2}\left( G_{R}\right)
}^{2}\right) +
\]
\[
+C_{1}e^{-\lambda \mu }\left\Vert Q\right\Vert _{H^{2}\left( G_{\mu
,R}\right) }^{2},\lambda \geq \lambda _{1}.
\]
Combining this with (\ref{26}), we obtain
\[
\int\limits_{G_{2\mu ,R}}\left[ \left( \Delta Q\right) ^{2}+\left( \nabla
Q\right) ^{2}+Q^{2}\right] e^{2\lambda r}dx\leq 
\]
\begin{equation}
\leq C_{1}\left( 1+\left\Vert V^{\ast }\right\Vert _{H^{2}\left(
G_{R}\right) }^{2}\right) \left( \delta ^{2}e^{2\lambda R}+e^{-\lambda \mu
}\right) .  \label{30}
\end{equation}
Choose $\lambda =\lambda \left( \delta \right) $ such that
\begin{equation}
\delta ^{2}e^{2\lambda \left( \delta \right) R}=\delta .  \label{300}
\end{equation}
Hence, 
\begin{equation}
\lambda \left( \delta \right) =\ln \left( \delta ^{-1/\left( 2R\right)
}\right) .  \label{31}
\end{equation}
5
We choose $\delta _{0}=\delta _{0}\left( L,G_{\mu ,R}\right) \in \left(
0,1\right) $ so small that $\ln \left( \delta _{0}^{-1/\left( 2R\right)
}\right) =\lambda _{1}.$ Hence by (\ref{31})
\begin{equation}
\lambda \left( \delta \right) >\lambda _{1},\forall \delta \in \left(
0,\delta _{0}\right) ,  \label{32}
\end{equation}
\begin{equation}
e^{-\lambda \left( \delta \right) \mu }=\delta ^{\mu /\left( 2R\right) }.
\label{33}
\end{equation}
Since $\mu \in \left( 0,R/2\right) ,$ then $\delta <\delta ^{\mu /\left(
2R\right) }.$ Hence, (\ref{30})-(\ref{33}) imply
\[
\left\Vert \Delta Q\right\Vert _{L^{2}\left( G_{2\mu ,R}\right) }+\left\Vert
Q\right\Vert _{H^{1}\left( G_{2\mu ,R}\right) }\leq C_{1}\left( 1+\left\Vert
V^{\ast }\right\Vert _{H^{2}\left( G_{R}\right) }\right) \delta ^{\mu
/\left( 4R\right) },\forall \delta \in \left( 0,\delta _{0}\right) .
\]%
This, (\ref{17}) and (\ref{18}) imply (\ref{160}). Remark 3.1 combined with (%
\ref{160}) obviously implies (\ref{161}). $\square $

\section{Numerical study}
\label{numer}

To compare the computational effectiveness of two methods mentioned above, we have conducted the series of numerical experiments. All computations have been carried out on the Lambda TensorBook with the 32 GB RAM, GPU NVIDIA RTX 2080 Max-Q with 8GB memory and 6-core Intel i7 running under Ubuntu 18.04. The graphical content of the paper have been produced with IDL 6.2.

\subsection{Data simulation}
\label{simul}

Let ${\it supp}(q(x))\subset\Omega\subset G_{R}\subset S$, where $C_{R}=\partial G_{R}$ and $C_{S}=\partial S$ are two 
concentric circles with the origin in ${\it supp}(q(x))$ and with the radii $R$ and $R_{S}$. Suppose that the point-like uniformly distributed transducers are situated on $C_{R}$. In acoustics the transducers are usually used as generation as reception of the acoustic waves, so that the transducer situated at $x_{0}\in C_{R}$ emits instantly an acoustic pulse while all other transducers situated at $x\in C_{R}, x\neq x_{0}$ receive the acoustic field. Such a geometrical configuration is typical in medical tomography and non-destructive testing.    

To simulate the data within the first model, we apply the finite element method, namely the Gal;erkin method, to the numerical solution of (\ref{e_1})-(\ref{e_3}). Let $\tau^{N}$ be a triangulation of $\Omega$ with $N$ nodes $x_{1},...,x_{k},...,x_{N}$ and $N_{tr}$ triangles $\Delta_{k}$. According to \cite{pest}, let $\{\psi_{n}(x)\}_{n=1}^{N}$ be a set of piecewise linear basis functions, such that $\psi_{n}(x_{m}) = \delta_{nm}, n,m=1,...,N$, and the solution $u^{f}$ of 
(\ref{e_1})-(\ref{e_3}) is approximated by the finite series
\begin{equation}
u^{f}_{N}(x,t) = \sum_{n=1}^{N}U^{f}_{n}(t)\psi_{n}(x).
\label{proj}
\end{equation}
In doing so, we approximate the coefficient $q$ by its discrete analogue with the value $q_{k}= const$ in the $k$-th triangle. Since the unknown coefficients $U^{f}_{n}(t)$  satisfy to a system of the second-order linear ODEs with the constant coefficients and homogeneous initial conditions, this system is solved by the standard methods. Note that the coefficients depend only on $q_{k}$ and $\psi_{n}$. The vector $U^{f}(2T)=(U^{f}_{1},..U^{f}_N))$ at the boundary nodes of $\tau^{N}$ represents the data used in the reconstruction procedure. The discrete Ricker wavelets and $C^{\infty}$-approximations of the Dirac delta functions are used as the controls $f$ at the boundary nodes. 

To simulate the data $g_{0}, g_{1}$ within the second model (\ref{c_1})-(\ref{c_2}) (or (\ref{c_11})-(\ref{c_21})), one may use at least two approaches. In the first approach, one applies the Laplace transform to the model, reduces it to the Lippmann-Schwinger equation, and solves numerically two coupled equations for each parameter of the Laplace transform in order to obtain its boundary values. 
However, in the course of numerical experiments, it was found that the accuracy of the numerical solution of the Lippmann-Schwinger equation turned out to be lower compared to the accuracy of the numerical solution to the Cauchy problem 
(\ref{c_1})-(\ref{c_2}) (or (\ref{c_11})-(\ref{c_21})). 
Therefore, we focus at the numerical solution to the latter. Clearly, since the use of the finite-difference methods requires a bounded computational domain, its boundary needs to be transparent to outgoing waves. To ensure this, the so-called absolute transparency boundary condition (see \cite{ sofron}) is posed on $C_{S}$ and the following initial boundary value problem is introduced instead of (\ref{c_1})-(\ref{c_2}) 
\begin{eqnarray}
& & q(x)v_{tt} - \Delta{v} = 0~x\in S, t > 0
\label{abc_1} \\
& & v(x,0) =0,~~ v_{t}(x,0) = \delta(x-x_{0}),~x\in S,x_{0}\in C_{R}, \label{abc_2}  \\
& & {\cal B}v = 0~\mbox{on}~C_{S}. \label{abc_3}
\end{eqnarray}
Here, ${\cal B}$ is a linear operator, such that for a certain $\hat{q} = const > 0$ resulted from the homogenization of $q(x)$in $S$ (e.g., the mean of $q(x)$), $u\simeq v$ in $S$. The algorithms for the construction of such operators ${\cal B}$ can be found in \cite{grot,sofron}. In computations, the Dirac delta function $\delta(u)$ is approximated by the $\varepsilon$-parametric "cap" functions (see, e.g., \cite{vlad}, p.8) with a sufficiently small parameter $\varepsilon > 0$. As previous, the problem (\ref{abc_1})-(\ref{abc_3}) is numerically solved on the disk $S$ for a sufficiently large $T$ by the Galerkin method. In the numerical experiments $T$ is chosen to provide $\|v(\cdot,\cdot;T)\|_{2} \le 10^{-3}$ in $S$. 
The numerical solution of this problem is used to compute the boundary conditions (\ref{bc_V}) on the sides of a square $D$ inscribed in $C_{R}$, so that the coefficient $q(x)$ is reconstructed in this square. By analogy with (\ref{noi}), a noise model is used to simulate the perturbed boundary data.

It follows from (\ref{v_f}) and (\ref{rel_4}) that the data $\psi(x,x_{0})$ in (\ref{lavren}) is determined by the Laplace transform $\tilde{u}(\cdot,\cdot;p)$. In computations the discrete truncated Laplace transform is evaluated, i.e., for each fixed $p\in(0,e^{-\gamma})$ a quadrature formula is used to compute the values of $\exp{(-pt)}u(\cdot,\cdot;p)$ on the segment $[\tau, T]$. To compute the discrete truncated Laplace transform, we use the fast algorithm by Rokhlin \cite{rokhl}, which requires $(n+m){\it log}_{2}(1/\epsilon_{pr})$ arithmetic operations, where $n,m$ are the dimensions in $t$ and $p$, and 
$\epsilon_{pr}$ is the precision of computations.       

\subsection{Inversion}
\label{discr}

\subsubsection{A discrete analogue of the BCM.}
\label{discr_bcm}

The discrete analogue of Inverse Problem 1 is formulated as: Given the discrete analogs 
$[f_{i},f_{j}]^{N_{b}}, i,j=1,...,N_{b}+1$ of the bilinear form
\[
[f,g] = \int_{\Omega}q(x)u^{f}(x,T)u^{g}(x,T)dx
\]
that is related to the Neumann-to-Dirichlet map, determine $q_{k}$ in $\tau^{N}$. Here, $N_{b}$ is the number of the boundary nodes. The detailed description of the inversion procedure can be found in \cite{pest}. According to the Pestov's approach, we obtain a system of linear algebraic equation
\begin{equation}
\sum_{k=1}^{N_{tr}}\left( \sum_{n=1}^{N}\sum_{m=1}^{N}a_{k}^{(n,m)}\varphi_{\alpha}^{(n)}\varphi_{\beta}^{(m)} \right) q_{k} = b_{\alpha\beta}
\label{syst}
\end{equation}
where
\[
a_{k}^{(n,m)} = \int_{\Delta_{k}}\psi_{n}(x)\psi_{m}(x)dx,~~b = [f_{\alpha},f_{\beta}],~\alpha, \beta = 1,...,N_{b}.
\]
The matrix in (\ref{syst}) is ill-conditioned. In our numerical experiments the condition number of the matrix in (\ref{syst}) grows dramatically with increasing $N$, attaining $10^{8}$ quickly for $N=64$. Therefore, we solve this system by minimizing the Tikhonov functional
\[
T[\hat{\rho}] = \frac{1}{2}\|A\hat{\rho}-b\|_{L^{2}}^{2} + \gamma\|\hat{\rho}\|_{H^{2}}^{2},~\gamma > 0.
\]
The standard methods (see, e.g., \cite{BG}) are used to solve numerically this problem.
\begin{figure}[!htb]
\centering
\includegraphics[scale=.75]{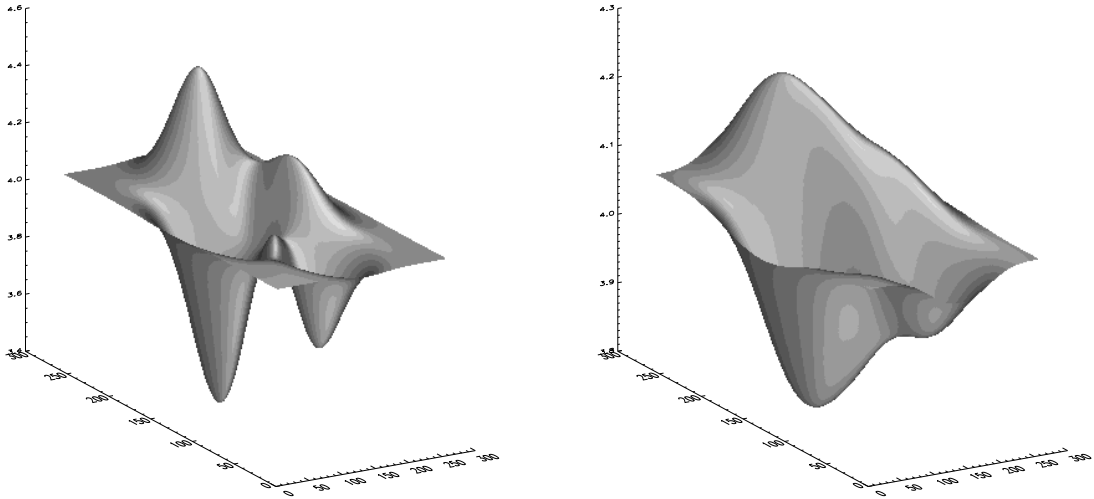}
\caption{Comparison of the shaded surfaces of $q$ reconstructed from the noiseless data. From left to right: the ground truth 
$q$, the reconstructions by the QRM and BCM. }
\label{figure_2}
\end{figure}

\subsubsection{A discrete approximation of the QRM.}
\label{imp_qrm}

We fix $x_{0}\partial G_{R}$ and approximate the linear differential operator (\ref{syst_V}) in a uniform grid of $G_{R}$ with a  standard second order 
finite difference analogue (see, e.g., \cite{sam}). 
Let $M$ be the resulting matrix. We also approximate the boundary conditions (\ref{bc_V}) with the second order central finite differences. The values of the normal derivative at the boundary nodes are computed through the internal grid nodes with the second order of accuracy. By adding the rows, containing numerical coefficients in the difference formulae of the boundary conditions, to the matrix $M$, we obtain a new matrix $A$, such that the central block of $A$ is the matrix $M$. As a result, a vector of unknowns $V$ is defined on both the internal and boundary 
grid nodes. This construction allows for approximating the constrained problem of minimizing the functional (\ref{0}) on $B$ 
by the unconstrained problem of minimizing the Tikhonov functional 
\begin{equation}
T_{\alpha}(V) = \|AV-b\|^{2}_{2} + \alpha\|V\|^{2}_{2},~\alpha > 0
\label{tikh}
\end{equation}
where the vector $b$ contains the values of $S_{0}$ and $S_{1}$, as well as zeros which correspond to the internal grid nodes. 
It should be emphasized that by virtue of (\ref{decomp}), $A$ and $b$ are just some approximations of the actual quantities. Therefore, for a fixed $\alpha > 0$ a solution $V_{\alpha}$ to the variational problem ${\it argmin}\{T_{\alpha}(V)\}$ satisfying the equation   
\begin{equation}
A^{t}AV_{\alpha} + \alpha E = A^{t}b,
\label{eul}
\end{equation}
where $E$ is a unit matrix, approximates the unique normal pseudo-solution of the system $AV=b$, which is, in general, incompatible. 
This is called the $\alpha$-approximation as opposed to its regularized solution. 
\begin{figure}[!htb]
\centering
\includegraphics[scale=.58]{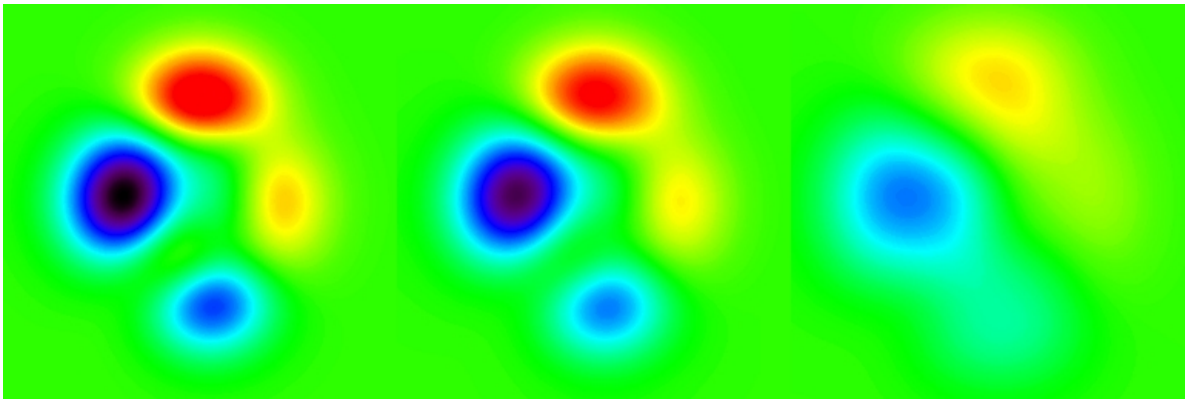}
\caption{Comparison of images of $q$ reconstructed from the noiseless data. From left to right: the ground truth 
$q$, the reconstructions by the QRM and BCM.}
\label{figure_3}
\end{figure}
The main uncertainty of this approach is that it is difficult to determine the acceptable error levels of $A$ and $b$ when selecting the regularization parameter $\alpha$ by a generalized residual principle (see, e.g., \cite{book_tgsy}, Ch.1, 2). On the other hand, the SVD-based minimal pseudo-inverse matrix method \cite{leo} has not provided the necessary robustness in our numerical experiments. For these reasons, we use a heuristic technique in order to select a quasi-optimal solution to the equation (\ref{eul}). Namely, starting with the sufficiently small $\alpha$, we solve (\ref{eul}) for $V_{\alpha}$ and then compute the approximate $\xi$ from a finite-difference analogue of (\ref{diff}) and (\ref{decomp}). Clearly, in this case for the noisy data we obtain a noisy image of $\xi$. To filter out the noise in an image, we have utilized an algorithm from image processing based on the sigma probability of a Gaussian distribution \cite{lee}. The main advantage of this algorithm is that it only exploits the data and does not need knowledge of the error levels. Also, it filters out the noise preserving image edges and fine details. In image processing this algorithm is deemed as the most computationally efficient noise filter. Gradually increasing the value of $\alpha$ and estimating the reconstruction errors in comparison with the ground truth models, we may select the parameter $\alpha$ that provides a quasi-optimal reconstruction of the coefficient $q$. 
Since the right-hand side of (\ref{lavren}) depends on $x_{0}$ as a parameter, the sought $\xi$ is found for each 
$x_{0}\in C_{R}$. After all such $\xi$'s are determined, one may either choose the best in some way, or calculate the average. In the numerical experiments we have found that the latter is the most suitable reconstruction.
%
\begin{figure}[!htb]
\centering
\includegraphics[scale=.58]{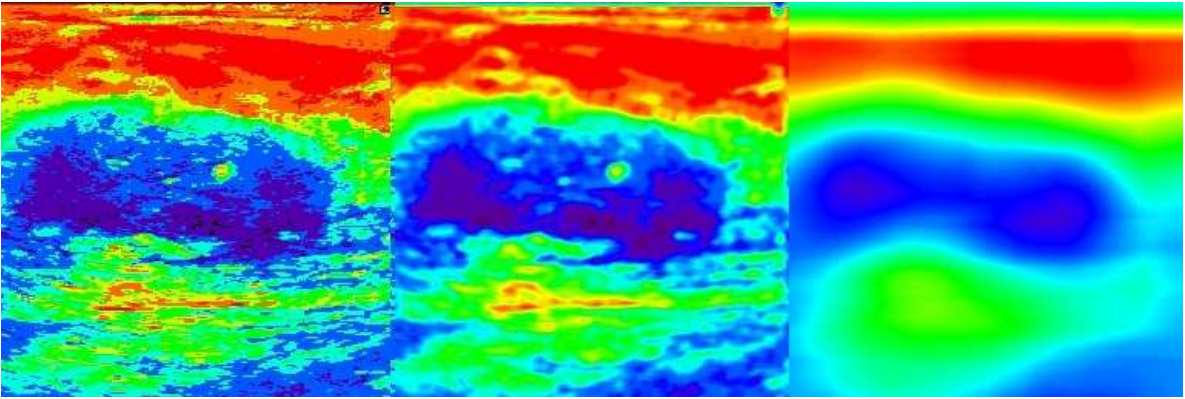}
\caption{Comparison of images of $q$ reconstructed from the noiseless data. From left to right: the discontinuous ground truth $q$; the reconstructions by the QRM ane BCM. The ground truth image represents a real ultrasound scan of breast borrowed from \cite{ultra}.}
\label{figure_6}
\end{figure}
            
\subsection{Some results of numerical experiments}
\label{numres}

In computations, we have used the dimensionless forward and inverse models, so that $R=1$ and the radius of $C_{S}$ was equal to 1.25.
Also, we note that the discrete analogue of $\{\Psi_{k}(x_{0})\}_{k=0}^{N}$ we use the SVD of a corresponding matrix. This ensures high accuracy of calculations of elements of the basis and, ultimately, robustness of reconstructions of the coefficient $q$.
As the model ground truth solutions $q$ of the inverse problems the smooth and discontinuous functions have been utilized. They model the distributions of $q$ in the fatty tissues, breast parenchyma and malignant breast lesions, the values of which 
are in the interval $(3.43, 4.53)$ of the dimensionless coefficient $q$. 
All reconstructions have been carried out on a uniform grid in $D$ with the number of nodes $256$ on each side of the square.    

\subsubsection{Reconstructions of a smooth distribution.}
\label{smoo}

%
\begin{figure}[!htb]
\centering
\includegraphics[scale=.75]{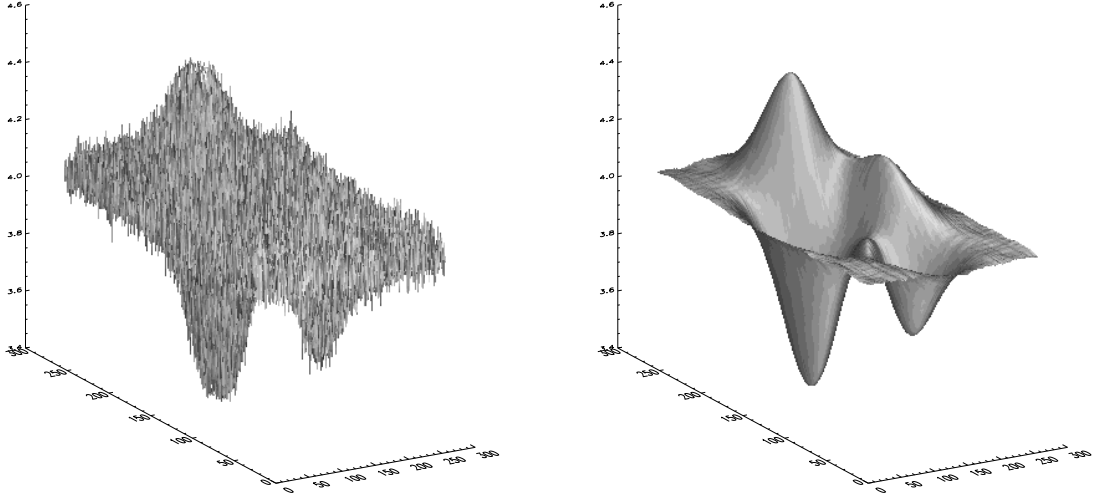}
\caption{The shaded surfaces of $q$ recovered by the QRM from the noisy data with the $1\%$ noise level. From left to right: the $\alpha$-approximation, $\alpha = 5\cdot 10^{-6}$; it is after the sigma-filtration.}
\label{figure_4}
\end{figure}
As a smooth ground truth solution, we use the following function
\begin{equation}
q(s,t) = c_{0} + c_{1}\left[ q_{1}(s,t) - q_{2}(s,t) - q_{3}(s.t) \right],~s,t\in[0,1],
\label{smo}
\end{equation}
where 
\begin{eqnarray*}
& & q_{1}(s,t) = c_{2}(1-3s)^{2}\exp{[-9s^{2}-(3t-2)^{2}]}, \\
& & q_{2}(s,t) = (\frac{3}{5}s-27s^{3}-(3(t-1))^{5})\exp{[-(9s^{2}+9(t-1)^{2})]}, \\
& & q_{3}(s,t) = \exp{[-(3s+1)^{2} - 9(t-1)^{2}]},
\end{eqnarray*}
and the positive constants $c_{0}, c_{1}, c_{2}$ are chosen so that $q(x_{1},x_{2}), (x_{1},x_{2})\in D$ possesses the mollifying property in a small neighborhood of $\partial D$. Here, $x_{1} = \sqrt{2}s - \sqrt{2}/2, x_{2}=\sqrt{2}t - \sqrt{2}/2$. The shaded surface of the ground truth $q$ is shown in Figure~\ref{figure_2} on the left. The shaded surface of $q$ reconstructed by the QRM from the noiseless data is depicted in Figure~\ref{figure_2} in comparison with the corresponding reconstruction by the BCM. Comparison of corresponding images is demonstrated in Figure~\ref{figure_3}. 
\begin{figure}[!htb]
\centering
\includegraphics[scale=.75]{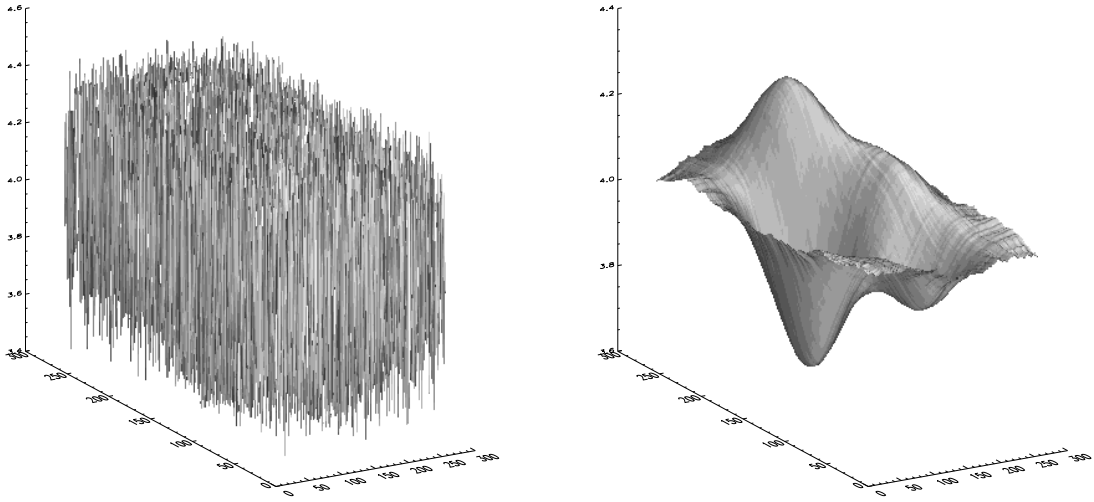}
\caption{The shaded surfaces of $q$ recovered by the QRM from the noisy data with the $5\%$ noise level. From left to right: its $\alpha$-approximation , $\alpha = 5\cdot 10^{-5}$; it is after the sigma-filtration.}
\label{figure_5}
\end{figure}

\subsubsection{Reconstructions of a discontinuous distribution.}
\label{disc}

As a discontinuous ground truth model, we utilize a real breast cancer image obtained by ultrasound scan (see \cite{ultra}). Since any digital image has the pixel structure, it can be represented by a discontinuous function. The image has been rescaled to the interval $(3.43, 4.53)$ that is the range of the dimensionless coefficient $q$. Comparison of reconstructions by the QRM and BCM is shown in Figure~\ref{figure_6}. 
   
\subsubsection{Robustness.}
\label{robs}

\begin{figure}[!htb]
\centering
\includegraphics[scale=.58]{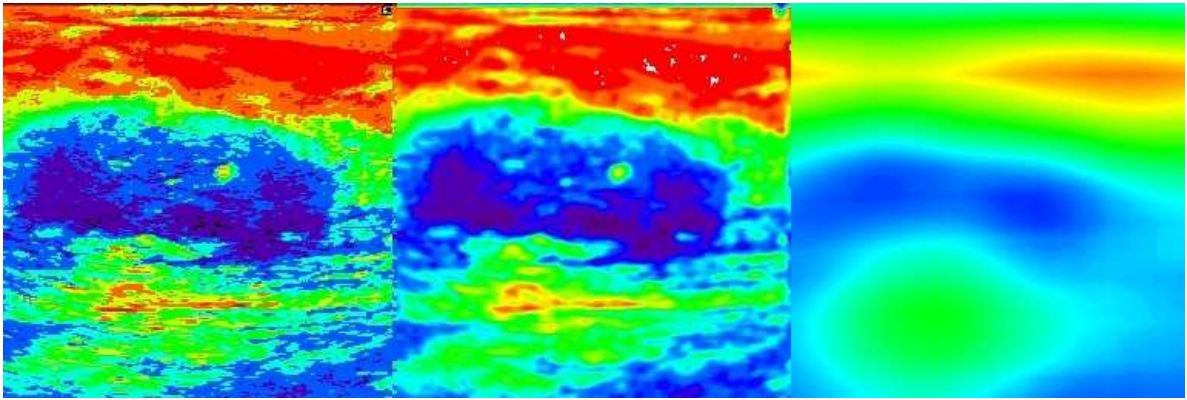}
\caption{Comparison of images of the discontinuous $q$ recovered from the noisy data with the $1\%$ noise level. From left to right: the ground truth $q$; the sigma-filtered reconstructions by the QRM and BCM. The ground truth image represents a real ultrasound scan of breast borrowed from \cite{ultra}.}
\label{figure_7}
\end{figure}
To simulate the noisy data within the first model, we use a stochastic model of the additive normally distributed noise
\begin{equation}
\tilde{U}^{f} = {U}^{f} + \delta\cdot\frac{\|{U}^{f}\|_{2}}{\|{\cal N}\|_{2}}\cdot{\cal N},
\label{noi}
\end{equation}
where $\delta > 0$ is a prescribed level of noise and ${\cal N}$ is the normally distributed pseudo-random vector with the zero mean and standard deviation 1. 

By analogy with (\ref{noi}) we simulate the noisy data (\ref{ch_1})-(\ref{ch_2}) for each 
$x_{0}\in C_{R}$ within the second model Due to the stochastic nature of (\ref{noi}) and the perturbed data (\ref{ch_1})-(\ref{ch_2}), the stochastic vectors have been generated 20 times to form a sample of reconstructed coefficients for each $\delta$. Then the results are represented by a mean over this sample. They are shown in Figures~\ref{figure_4} and \ref{figure_5} for the different levels of noise.   
\begin{figure}[!htb]
\centering
\includegraphics[scale=.58]{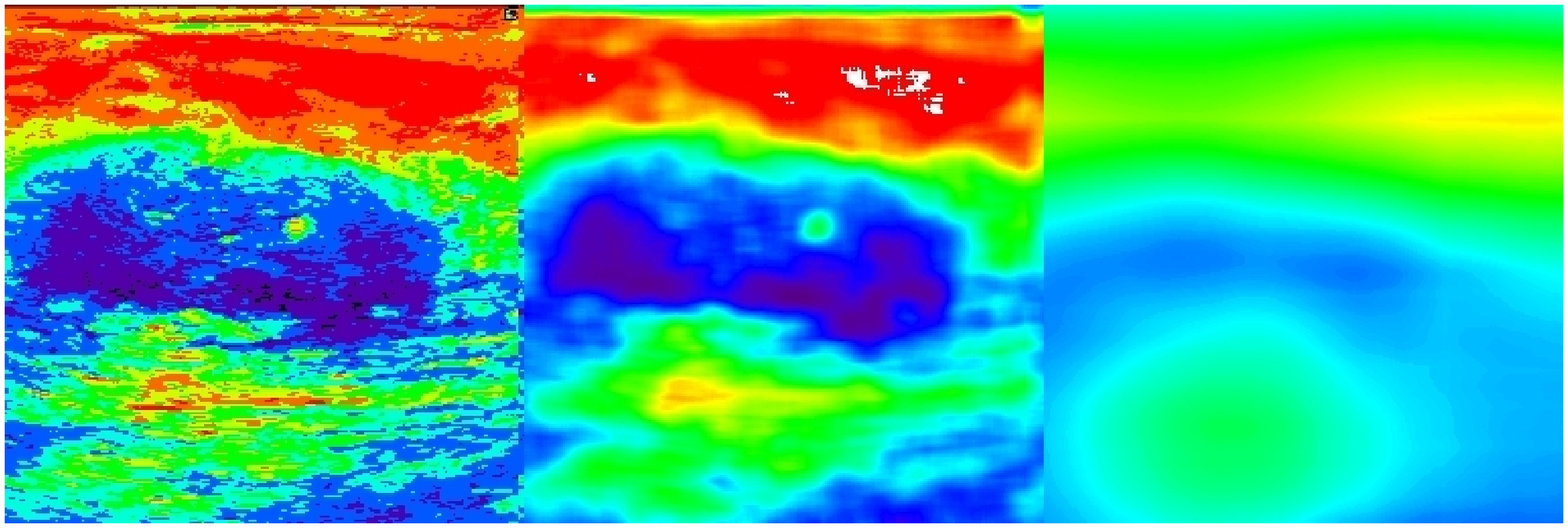}
\caption{ Comparison of images of the discontinuous $q$ recovered from the noisy data with the $5\%$ noise level. From left to right: the ground truth $q$; the sigma-filtered reconstructions by the QRM and BCM. The ground truth image represents a real ultrasound scan of breast borrowed from \cite{ultra}. }
\label{figure_8}
\end{figure}   
In the case of the discontinuous ground truth, the results of recovering images of $q$ from the noisy data are shown in Figures~\ref{figure_7} and \ref{figure_8}.      

\section{Conclusions}
\label{concl}

We carried out the comparative study of the QRM and BCM in two dimensions. The numerical implementations of these method were developed and used to demonstrate their computational effectiveness in the numerical experiments with the smooth and discontinuous coefficients of a wave equations. Under comparable conditions, which are the delta-like sources for the QRM and delta-like controls for the BCM, it was demonstrated in the numerical experiments that from a computational point of view, the QRM based on the Lavrentiev equation is superior to the BCM.  

\noindent
{\bf References}
\\

\end{document}